\documentclass[12pt,leqno]{article}
\usepackage{amsmath,amssymb,bbm,array}
\usepackage{amsthm}
\usepackage[all,cmtip]{xy}

\def\vds{\hbox{\hbox to 2\arraycolsep{\hss\vbox{\vbox to 1.3ex{ 
\vss\hbox{.}\vspace{.33ex}\hbox{.}\vspace{.33ex}\hbox{.}\vspace{.33ex}\hbox{.}\vspace{.33ex}\hbox{.}\vss}}\hss}}}

\makeatletter
\renewcommand{\section}{\@startsection{section}{1}{0mm}{5mm}{3mm}{\raggedright\bf\large}}

\def\@citex[#1]#2{\if@filesw\immediate\write\@auxout{\string\citation{#2}}\fi
  \def\@citea{}\@cite{\@for\@citeb:=#2\do
    {\@citea\def\@citea{\@citesep}\@ifundefined
       {b@\@citeb}{{\bf ?}\@warning
       {Citation `\@citeb' on page \thepage \space undefined}}%
{\csname b@\@citeb\endcsname}}}{#1}}
\def\@citesep{; }
\makeatother

\newtheoremstyle{Kang}{}{}{\itshape}{}{\bf}{}{.5em}{}
\theoremstyle{Kang}
\newtheorem{theorem}{Theorem}[section]

\newtheorem{corollary}[theorem]{Corollary}

\newtheoremstyle{Kremark}{}{}{}{}{\bf}{}{.5em}{}
\theoremstyle{Kremark}

\newtheorem*{remark}{Remark.}
\newtheorem{other}{}

\allowdisplaybreaks[1]  
\numberwithin{equation}{section}

\begin{document}

\title{Rationality for subgroups of $S_6$}
\author{Jian Zhou \\[3mm]
School of Mathematical Sciences, Peking University, Beijing,\\
E-mail: zhjn@math.pku.edu.cn }
\date{}

\maketitle


\setlength{\parskip}{2pt}

\begin{abstract}
{\bf Abstract.} For a transitive subgroup $G \le S_6$ which contain
$C_3 \times C_3$ as subgroup, we prove that $K(x_1,\dots,x_6)^G$ is
rational over $K$, where $K$ is any field, and $G$ acts naturally on
$K(x_1,\dots,x_6)$ by permutations on the variables. We also give an
application on construction of generic polynomials.
\end{abstract}

{\small \hspace*{4mm}2010 Mathematics Subject Classification.
Primary 12F20, 13A50, 14E08.

\hspace*{4mm}{\it Keywords}: rationality, symmetric group, Noether's problem, generic polynomial.\\
}

\section{Introduction}

Let $K$ be a field, and $x_1,\dots,x_n$ be variables. We have a
natural action of symmetric group $S_n$ on these $x_i$'s by
permutations, which can be extended to a natural action of $S_n$ on
the rational function field $K(x_1,\dots,x_n)$, with trivial action
of $G$ on $K$. For a subgroup $G \le S_n$, Noether asked if the
subfield $K(x_1,\dots,x_n)^G$ is rational (= purely transcendental)
over $K$, this problem has its origin in the study of inverse Galois
problem, and has been investigated by many people from different
aspects, see e.g. \cite{Sw} for a survey of related topics.

Our main result in this paper is

\begin{theorem} Let $K$ be any field, $G$ be a transitive subgroup
of $S_6$ which contains a subgroup isomorphic to $C_3 \times C_3$.
If $G \ne A_6$, then $K(x_1,\dots,x_6)^G$ is rational over $K$.
\end{theorem}

Recently, M. Kang and B. Wang \cite{KW} has investigated this
retionality problem for $S_n$ with $n \le 5$ and $n=7$, and get
fairly complete results.

For an arbitrary finite group $G$, consider the regular embedding of
$G$ into $S_{|G|}$ (where $|G|=$ number of elements in $G$), namely
$G$ acting on itself by left multiplication. This make $G$ a
transitive subgroup of $S_{|G|}$, and the commonly called
``Noether's problem" is the following:

{\it Let $\{x_g\}_{g \in G}$ be a family of variables, indexed by
$G$, and let $G$ act on these $x_g$'s by $h(x_g) = x_{hg},\,
\forall\, h,g \in G$. For a given field $K$, consider the rational
function field $K(x_g: g \in G)$ with natural $G$-action as above.
When does the fixed subfield $K(G)= K(x_g: g\in G)^G$ again rational
over $K$ ?}

Results on Noether's problem are abound, we just list \cite{Len},
\cite{Sa2}, \cite{HuK} for a glance.

From Theorem 1.1, we deduce

\begin{corollary}
Let $G$ be a transitive subgroup of $S_6$ which contains $C_3 \times
C_3$, and assume $G \ne A_6, S_6$, then Noether's problem has a
positive answer over any field $K$ of characteristic $\,\ne 2,3$.
For $G = S_6$, the same result holds provided that ${\rm char}(K)
\ne 2,3,5$.
\end{corollary}

The proof is standard, namely, we can embed the natural
$G$-representation $V_0 = \bigoplus_{1 \le i \le 6} K x_i$ into the
regular $G$-representation $V_{\rm reg} = \bigoplus_{g \in G} K x_g$
(the existence of such an embedding is explained at the end of \S2),
and appeal to the following basic tool (take $L = K(x_1,\dots,x_6)
\subseteq K(x_g: g \in G)$ under the above embedding):

\begin{theorem}
{\rm (\cite{HaK}, Th. 1)} Let $G$ be a finite group acting on
$L(y_1,\dots,y_m)$, the rational function field over a field $L$
with variables $y_1,\dots,y_m$. Assume that $L$ is stable under $G$,
and the restricted action of $G$ on $L$ is faithful. Moreover,
assume for any $\sigma \in G$, $$
\begin{pmatrix}
\sigma(y_1) \\ \vdots \\ \sigma(y_m)
\end{pmatrix} = A(\sigma) \begin{pmatrix}
y_1 \\ \vdots \\ y_m
\end{pmatrix} + B(\sigma), $$ where $A(\sigma) \in {\rm GL}_m(L)$ and
$B(\sigma) \in L^m$. Then $L(y_1,\dots,y_m) = L(z_1,\dots,z_m)$ for
some $z_1,\dots,z_m$ which are fixed by $G$. In particular,
$L(y_1,\dots,y_m)^G = L^G(z_1,\dots,z_m)$ is rational over $L^G$.
\end{theorem}

In \S2, we will describe the structure of those subgroups appeared
in Theorem 1.1, then in \S3, we present the proof of the theorem by
explicitly write down the generators for the invariants. In \S4, we
show that the Galois descent gives another approach, and finally in
\S5, we discuss the construction of generic polynomials for these
subgroups of $S_6$, as an application of our main theorem.

\section{Subgroups of $S_6$}

Transitive subgroups of $S_6$ (other than $A_6$ and $S_6$) are
listed in \cite{DM}, page 60. Up to conjugacy, there are 14 of them,
numbered by T6.1-T6.14 (refer to \cite{DM} for the descriptions of
all these groups by generators). We can divide the transitive
subgroups of $S_6$ into 3 types,

Type \; I : \ T6.1 -- T6.8;

Type \;II : \ T6.9 -- T6.12;

Type III : \ T6.13, T6.14, $A_6$, and $S_6$.

Among these, Type I and II are solvable groups, while Type II
consists of those with a $C_3 \times C_3$ as subgroup. In fact,
groups T6.1 to T6.8 are of orders 6,6,12,48,24,24,12,24, so they
won't contain $C_3 \times C_3$ by a trivial reason.

In this paper, we will consider the rationality problem for Type II
groups only. The rationality problem for Type I subgroups may be
investigated from a more general setting for wreath products and
will be discussed in a separate paper.

Type III are unsolvable groups. Except for $S_6$ itself, the
rationality problem becomes very difficult, but see a related result
of Maeda \cite{Mae}.

Let us now give a detailed description of Type II groups:

Firstly, in $S_6$ we let $$ \sigma_1 = (123), \quad \sigma_2 =
(456), \quad \tau = (14)(25)(36), \quad \lambda = (1425)(36), $$
note that $\sigma_2 = \tau \sigma_1 \tau = \lambda \sigma_1
\lambda^{-1}$ and $\lambda^2 = (12)(45)$. Now $\sigma_1, \sigma_2$
generate a subgroup $C_3 \times C_3 < S_6$, which will be contained
in each of our Type II groups.

$\bullet$ T6.9, referred to $G_1$ for simplicity. This is a group of
order 72, with generators $\sigma_1, \tau$, and $(12)$. Since $\tau
(12) \tau = (45)$, and $\sigma_1, (12)$ (resp. $\sigma_2, (45)$)
generate $S_3$ on $1,2,3$ (resp. on $4,5,6$), we see easily that
$S_3 \times S_3$ is a normal subgroup of $G_1$ of order 36, and
$G_1/(S_3 \times S_3)$ is generated by $\tau$.

$\bullet$ T6.10, referred to $G_2$. This is a group of order 36,
with generators $\sigma_1$ and $\lambda = (12)\tau$, and it is the
intersection of $G_1$ with $A_6$.\footnote{We remark that there is a
misprint in \cite{DM} for T6.10 (page 60), the $\lambda' =
(1542)(36)$ there should be replaced by $\lambda = (1425)(36)$ here,
for otherwise $\lambda' \sigma_1 \lambda' = (14652)$ is of order 5
in T6.10, absurd.} One checks that $C_3 \times C_3$ generated by
$\sigma_1,\sigma_2$ is a normal subgroup of $G_2$ of order 9, and
$G_2/(C_3 \times C_3) \simeq C_4$ is generated by $\lambda$.

$\bullet$ T6.11, referred to $G_3$. This is a group of order 36,
with generators $\sigma_1, \tau, \lambda^2$. One checks similarly
that $C_3 \times C_3$ is a normal subgroup of $G_3$, and $G_3/(C_3
\times C_3) \simeq C_2 \times C_2$ is generated by $\tau,
\lambda^2$.

$\bullet$ T6.12, referred to $G_4$. This is a group of order 18,
with generators $\sigma_1, \tau$. Therefore $G_4/(C_3 \times C_3)
\simeq C_2$ is generated by $\tau$.

\begin{remark} For $G = G_1,G_2,G_3$ or $G_4$, and for a field $K$
of characteristic $\ne 2,3$, the regular $K$-representation $V_{\rm
reg}$ of $G$ is semisimple, and contains all the irreducible one's,
so to prove that the natural representation $V_0 = \bigoplus_{1 \le
i \le 6} K x_i$ is embedded in $V_{\rm reg}$, it suffices to show
that all the irreducible components of $V_0$ are of multiplicity
one.

We may assume that all the $3$-th roots of unity are in $K$, then it
is easily seen that $V_0$ as a representation of $C_3 \times C_3 =
\langle \sigma_1,\sigma_2 \rangle$ is decomposed into six
1-dimensional subrepresentations, two of them are trivial
representations, while the other four are nontrivial and mutually
non-isomorphic. Now the two trivial ones are interchanged by $\tau$
(resp. $\lambda$), so their direct sum as a representation of $G_4$
(resp. $G_2$) can be decomposed into a trivial representation and a
nontrivial one, both of dimension 1. From these we see easily that
the multiplicity one property holds for $V_0$ as a representation of
$G_4 < G_3 < G_1$ or $G_3$.

The case $G = S_6$ is similar, but we need ${\rm char}(K) \ne 2,3,5$
to ensure the semisimplicity of $V_{\rm reg}$.
\end{remark}

\section{Proof of Theorem 1.1}

We will prove our main theorem through a case by case analysis. For
$S_6$ itself, it is a classical result, so we only need to consider
the rationalities for $G_1,G_2,G_3,G_4$.

As is expected, we will make use of Masuda's formula for invariants
of $C_3$:

\begin{theorem} {\rm (\cite{Mas}, Th. 3, \cite{HoK}, Th. 2.2)}
Let $K$ be any field, $x,y,z$ be variables. Consider the natural
action of $S_3$ on $K(x,y,z)$ by permutation of variables, and let
$A_3 (\simeq C_3) < S_3$ be the alternating subgroup, then
$K(x,y,z)^{A_3}$ is rational over $K$, and equals to $K(x+y+z,u,v)$,
where
$$
u = \frac{x^2y+y^2z+z^2x-3xyz}{x^2+y^2+z^2-xy-yz-zx}, \quad v =
\frac{xy^2+yz^2+zx^2-3xyz}{x^2+y^2+z^2-xy-yz-zx}.
$$
\end{theorem}

\bigskip

{\bf Rationality for $G_1$}

Let $y_1,y_2,y_3$ (resp. $y_4,y_5,y_6$) be the 3 elementary
polynomials for $x_1,x_2,x_3$ (resp. $x_4,x_5,x_6$), then $$
K(x_1,x_2,x_3,x_4,x_5,x_6)^{S_3 \times S_3} =
K(y_1,y_2,y_3,y_4,y_5,y_6),
$$
and $\tau$ acts on the $y_i$'s as follows: $$ \tau : \; y_1
\leftrightarrow y_4,\; y_2 \leftrightarrow y_5,\; y_3
\leftrightarrow y_6. $$ Let $$ \begin{array}{lll} z_1 = y_1+y_4,
\;\;\; & z_2 = y_2+y_5, \;\;\; & z_3 = y_3+y_6, \\ z_4 = y_1y_4,
\;\;\; & z_5 = y_1y_5+y_4y_2, \;\;\; & z_6 = y_1y_6+y_4y_3,
\end{array} $$
one checks easily that $$
K(y_1,y_2,y_3,y_4,y_5,y_6)^{\langle\tau\rangle} =
K(z_1,z_2,z_3,z_4,z_5,z_6). $$ Note that this is valid for any
fields, including those of characteristic 2, therefore
$K(x_1,x_2,x_3,x_4,x_5,x_6)^{G_1} = K(z_1,z_2,z_3,z_4,z_5,z_6)$ is
rational over $K$.

\begin{remark}
Similar method applies to the subgroup $H < S_{2n}$ generated by
$(12)$, $(12\dots n)$, and $(1,n+1)(2,n+2)\cdots(n,2n)$. In fact,
$H$ has a normal subgroup $S_n \times S_n$, and $H/(S_n \times S_n)$
is generated by $(1,n+1)(2,n+2)\cdots(n,2n)$.
\end{remark}

{\bf Rationality for $G_4$ and $G_3$}

In this case, we let \begin{equation} \tag{$*$} \left\{
\begin{array}{ll} u_1 = x_1+x_2+x_3, \; \\ u_2 =
\frac{x_1^2x_2+x_2^2x_3+x_3^2x_1-3x_1x_2x_3}{x_1^2+x_2^2+x_3^2-x_1x_2-x_2x_3-x_3x_1},
\; & u_3 =
\frac{x_1x_2^2+x_2x_3^2+x_3x_1^2-3x_1x_2x_3}{x_1^2+x_2^2+x_3^2-x_1x_2-x_2x_3-x_3x_1},
\\ u_4 = x_4+x_5+x_6, \; \\ u_5 =
\frac{x_4^2x_5+x_5^2x_6+x_6^2x_4-3x_4x_5x_6}{x_4^2+x_5^2+x_6^2-x_4x_5-x_5x_6-x_6x_4},
\; & u_6 =
\frac{x_4x_5^2+x_5x_6^2+x_6x_4^2-3x_4x_5x_6}{x_4^2+x_5^2+x_6^2-x_4x_5-x_5x_6-x_6x_4}.
\end{array} \right. \end{equation}
By Masuda's formula, we have $$ K(x_1,x_2,x_3,x_4,x_5,x_6)^{C_3
\times C_3} = K(u_1,u_2,u_3,u_4,u_5,u_6),
$$ and $\tau$ acts on the $u_i$'s as $u_1
\leftrightarrow u_4,\, u_2 \leftrightarrow u_5,\, u_3
\leftrightarrow u_6$. Therefore if we set $$ \begin{array}{lll} v_1
= u_1+u_4, \;\;\; & v_2 = u_2+u_5, \;\;\; & v_3 = u_3+u_6, \\ v_4 =
u_1u_4, \;\;\; & v_5 = u_1y_5+u_4y_2, \;\;\; & v_6 = u_1u_6+u_4u_3,
\end{array} $$ then $K(x_1,x_2,x_3,x_4,x_5,x_6)^{G_4} =
K(v_1,v_2,v_3,v_4,v_5,v_6)$.

$G_4$ is a normal subgroup of $G_3$, and $G_3/G_4$ is of order 2,
with generator $\lambda^2$. Therefore with notations as above, $$
K(x_1,x_2,x_3,x_4,x_5,x_6)^{G_3} =
K(v_1,v_2,v_3,v_4,v_5,v_6)^{\langle\lambda^2\rangle}.
$$ The action of $\lambda^2$ on $v_i$'s is given by $$ \lambda^2 : \; v_2
\leftrightarrow v_3,\; v_5 \leftrightarrow v_6,\; \text{and} \;v_1,
v_4 \; \text{fixed},
$$ so we see that $$ K(v_1,v_2,v_3,v_4,v_5,v_6)^{\langle\lambda^2\rangle} =
K(v_1,v_4,v_2+v_3,v_5+v_6,v_2v_3,v_2v_6+v_3v_5) $$ is rational over
$K$.

\bigskip

{\bf Rationality for $G_2$}

Since $C_3 \times C_3$ is also a normal subgroup of $G_2$, we know
as above $$ K(x_1,x_2,x_3,x_4,x_5,x_6)^{C_3 \times C_3} =
K(u_1,u_2,u_3,u_4,u_5,u_6), $$ where $u_1,u_2,u_3,u_4,u_5,u_6$ are
given by $(*)$. Now $G_2/(C_3 \times C_3)$ is generated by $\lambda
= (1425)(36)$, and $\lambda$ acts on $u_i$'s as follows:
$$ \lambda :\; u_1 \leftrightarrow u_4,\; u_2 \mapsto u_5 \mapsto
u_3 \mapsto u_6 \mapsto u_3. $$ In this stage, we need to treat the
following two cases separately, depending on ${\rm char}(K)$ equals
to 2 or not.

\bigskip

{\it Case} A --- ${\rm char}(K) \ne 2$.

\bigskip

Consider the fixed field of $\lambda^2 = (12)(45)$\, first.
$\lambda^2$ fixes $u_1,u_4$, and exchanges $u_2$ and $u_3$ (resp.
$u_5$ and $u_6$), so if we let
$$
\left\{ \begin{array}{lll} w_1 = u_1,
\;\;\; & w_2 = u_2+u_3, \;\;\; & w_3 = (u_2-u_3)/(u_5-u_6), \\
w_4 = u_4, \;\;\; & w_5 = u_5+u_6, \;\;\; & w_6 =
(u_2-u_3)(u_5-u_6),
\end{array} \right. $$ it is easily seen that $$
K(u_1,u_2,u_3,u_4,u_5,u_6)^{\langle\lambda^2\rangle} =
K(w_1,w_2,w_3,w_4,w_5,w_6)
$$ (here we need the assumption ${\rm char}(K) \ne 2$). Now $\lambda$ is
represented by $$ \lambda :\; w_1 \leftrightarrow w_4,\; w_2
\leftrightarrow w_5,\; w_3 \mapsto -\frac1{w_3},\; w_6 \mapsto -w_6.
$$ To get the final result, we put $$ \left\{
\begin{array}{llll} t_1 = w_1-w_4,
\;\; & t_2 = w_2-w_5, \;\; & t_3 = w_3-\frac{1}{w_3}, \;\; & s_3 = w_3+\frac{1}{w_3}, \\
t_4 = w_1+w_4, \;\; & t_5 = w_2+w_5, \;\;\; & t_6 = w_6,
\end{array} \right. $$ we see $K(w_1,w_2,w_3,w_4,w_5,w_6) = K(t_1,t_2,t_3,s_3,t_4,t_5,t_6)$,
and $\lambda(t_1) = -t_1$, $\lambda(t_2) = -t_2$, $\lambda(t_3) =
t_3$, $\lambda(s_3) = -s_3$, $\lambda(t_4) = t_4$, $\lambda(t_5) =
t_5$, $\lambda(t_6) = -t_6$. So $K(w_1,w_2,w_3,w_4,$
$w_5,w_6)^{\langle\lambda\rangle} =
K(s_3t_1,s_3t_2,t_3,t_4,t_5,s_3t_6)$ is rational over $K$.

\bigskip

{\it Case} B --- ${\rm char}(K) = 2$.

\bigskip

Here we let $$ \begin{array}{lll} w_1 = u_1,
\;\;\; & w_2 = u_2+u_3, \;\;\; & w_3 = u_2u_3, \\
w_4 = u_4, \;\;\; & w_5 = u_5+u_6, \;\;\; & w_6 = u_2u_6+u_3u_5,
\end{array} $$ and it gives $K(u_1,u_2,u_3,u_4,u_5,u_6)^{\langle\lambda^2\rangle} =
K(w_1,w_2,w_3,w_4,w_5,w_6)$. By some calculations, we find $$
\lambda :\; w_1 \leftrightarrow w_4,\;\; w_2 \leftrightarrow
w_5,\;\; w_6 \mapsto w_6 + w_2w_5,\;\; w_3 \mapsto \frac{w_5^2 w_3 +
w_6^2 + w_2w_5w_6}{w_2^2}. $$ To deal with this action, we make the
following change of variables:
$$ t_1 = w_1+w_4,\;\; t_2 = w_2+w_5,\;\; t_3 = \frac{w_5}{w_2}w_3,\;\;
t_4 = w_2w_1+w_5w_4,\;\; t_5 = w_5,\;\; t_6 = w_6. $$ $\lambda$\,
acts on the $t_i$'s as follows: $$ \begin{array}{ll} \lambda : & t_1
\mapsto t_1,\;\; t_4 \mapsto t_4,\;\; t_2 \mapsto t_2,\;\; t_5
\mapsto t_5 + t_2, \\ & \displaystyle t_6 \mapsto t_6 + t_2t_5 +
t_5^2,\;\; t_3 \mapsto t_3 + \frac{t_6(t_6 + t_2t_5 + t_5^2)}{t_2t_5
+ t_5^2}.
\end{array} $$ Note the special role played by $t_2$ and $t_5$,
we finally have \begin{multline*}
K(u_1,u_2,u_3,u_4,u_5,u_6)^{\langle\lambda\rangle} =
K(w_1,w_2,w_3,w_4,w_5,w_6)^{\langle\lambda\rangle} =
K(t_1,t_2,t_3,t_4,t_5,t_6)^{\langle\lambda\rangle} \\ =
K\bigg(t_1,\, t_4,\, t_2,\, t_2t_5+t_5^2,\, t_6 + (t_2t_5 + t_5^2)
\, \frac{t_5}{t_2},\, t_3 + \frac{t_6(t_6 + t_2t_5 + t_5^2)}{t_2t_5
+ t_5^2} \, \frac{t_5}{t_2}\bigg),
\end{multline*} which is rational over $K$.

\section{Another method}

The method used in the previous section relies heavily upon our
knowledge of $C_3$-invariants (i.e., Masuda's formula). There is
another method, i.e., Galois descent, which does not require such
explicit results on special invariants, and is applicable to more
general groups. See Swan's report \cite{Sw} for an introduction.

ÔÚ ±¾½ÚÖУ¬ we show that the descent method works equally well for
the groups $G_2,G_3,G_4$. For simplicity, let's consider only $G_3$,
to indicate the method. More precisely, we add a primitive $3$-th
root of unity $\zeta_3$ to the base field $K$ (assuming ${\rm
char}(K) \ne 3$), and consider the fixed subfield of
$K(\zeta_3)(x_1,\dots,x_6)$ under $G_3$ and the Galois group ${\rm
Gal}(K(\zeta_3)/K)$. For field $K$ of characteristic 3, we will also
give a new proof here.

Recall that $G_3$ is generated by $\sigma_1 = (123)$, $\tau =
(14)(25)(36)$, and $\lambda^2 = (12)(45)$, and it has a normal
subgroup $C_3 \times C_3$ generated by $\sigma_1$ and $\sigma_2 =
\tau\sigma_1\tau = (456)$. Without consulting Masuda's formula, we
need to put the action of $\sigma_1$ and $\sigma_2$ into a diagonal
form (or triangular form when ${\rm char}(K) = 3$).

\bigskip

{\it Case} A --- ${\rm char}(K) \ne 3$.

\bigskip

Let $$ \left\{ \begin{array}{l} y_1 = x_1 + x_2 + x_3, \\ y_2 =
\zeta_3^2 x_1 + \zeta_3 x_2 + x_3, \\ y_3 = \zeta_3 x_1 + \zeta_3^2
x_2 + x_3,
\end{array} \right. \quad \text{and} \quad
\left\{ \begin{array}{l} y_4 = x_4 + x_5 + x_6, \\ y_5 = \zeta_3^2
x_4 + \zeta_3 x_5 + x_6, \\ y_6 = \zeta_3 x_4 + \zeta_3^2 x_5 + x_6.
\end{array} \right. $$ One verifies that
$$ \begin{array}{ll} \sigma_1 : & y_2 \mapsto \zeta_3 y_2,\; y_3 \mapsto \zeta_3^2 y_3,\;
\text{and} \; y_1,y_4,y_5,y_6,\zeta_3 \;\text{fixed,} \\ \sigma_2 :
& y_5 \mapsto \zeta_3 y_5,\; y_6 \mapsto \zeta_3^2 y_6,\; \text{and}
\; y_1,y_2,y_3,y_4,\zeta_3 \;\text{fixed,} \\ \,\tau \;: & y_1
\leftrightarrow y_4,\; y_2 \leftrightarrow y_5,\; y_3
\leftrightarrow y_6,\; \text{and} \; \zeta_3 \;\text{fixed,} \\
\lambda^2 : & y_2 \leftrightarrow y_3,\; y_5 \leftrightarrow y_6,\;
\text{and} \; y_1,y_4,\zeta_3 \;\text{fixed.} \end{array}
$$ It is now easily checked
$$ K(\zeta_3)(y_1,\dots,y_6)^{C_3 \times C_3} =
K(\zeta_3)(z_1,\dots,z_6), $$ where $z_1 = y_1, z_2 = y_2^2/y_3, z_3
= y_3^2/y_2, z_4 = y_4, z_5 = y_5^2/y_6, z_6 = y_6^2/y_5$. And we
have
$$ \begin{array}{ll} \,\tau \;: & z_1 \leftrightarrow z_4,\; z_2 \leftrightarrow z_5,\;
z_3 \leftrightarrow z_6,\; \text{and} \; \zeta_3 \;\text{fixed,} \\
\lambda^2 : & z_2 \leftrightarrow z_3,\; z_5 \leftrightarrow z_6,\;
\text{and} \; z_1,z_4,\zeta_3 \;\text{fixed.}
\end{array} $$ Therefore $$ K(\zeta_3)(z_1,\dots,z_6)^{\langle \tau
\rangle} = K(\zeta_3)(u_1,\dots,u_6), $$ with $u_1 = z_1 + z_4, u_2
= z_2 + z_5, u_3 = z_3 + z_6, u_4 = z_1 z_4, u_5 = z_1 z_5 + z_4
z_2, u_6 = z_1 z_6 + z_4 z_3$, and $$ \lambda^2 :\; u_2
\leftrightarrow u_3,\; u_5 \leftrightarrow u_6,\; \text{and} \;
u_1,u_4,\zeta_3 \;\text{fixed.} \qquad\qquad $$

Now we come to the descent step. If $\zeta_3 \in K$, simply write
$$ K(\zeta_3)(u_1,\dots,u_6)^{\langle \lambda^2
\rangle} = K(u_1,u_4,u_2+u_3,u_5+u_6,u_2u_3,u_2u_5+u_3u_6), $$ we're
done. If $\zeta_3 \notin K$, then ${\rm Gal}(K(\zeta_3)/K) =
\{1,\rho\}$, with $\rho : \zeta_3 \mapsto \zeta_3^2$. Extend the
action of $\rho$ onto $K(\zeta_3)(x_1,\dots,x_6)$ by setting
$\rho(x_i) = x_i$, we find that $$ \rho \;:\; u_2 \leftrightarrow
u_3,\; u_5 \leftrightarrow u_6,\; \zeta_3 \mapsto \zeta_3^2,\;
\text{and} \; u_1,u_4 \;\text{fixed.} \quad $$ And the action of
$\rho$ and $G_3$ on $K(\zeta_3)(x_1,\dots,x_6)$ commutes. Therefore
$$
\begin{array}{rcl}
K(x_1,\dots,x_6)^{G_3} &=& K(\zeta_3)(x_1,\dots,x_6)^{G_3 \times
\langle\rho\rangle}
\\ &=& K(\zeta_3)(u_1,\dots,u_6)^{\langle \lambda^2 \rangle \times
\langle\rho\rangle } \\ &=& K(\zeta_3)(u_1,\dots,u_6)^{\langle
\lambda^2\rho \rangle \times \langle \lambda^2 \rangle} \\ &=&
K(u_1,\dots,u_6)^{\langle \lambda^2 \rangle} \qquad \text{(since
$\lambda^2\rho$ fixes $u_i$'s, and maps $\zeta_3$ to $\zeta_3^2$)} \\
&=& K(u_1,u_4,u_2+u_3,u_5+u_6,u_2u_3,u_2u_5+u_3u_6) .
\end{array} $$

{\it Case} B --- ${\rm char}(K) = 3$.

\bigskip

Letting $$ \left\{ \begin{array}{l} y_1 = x_1 + x_2 + x_3, \\ y_2 =
- x_1 + x_2, \\ y_3 = x_1,
\end{array} \right. \quad \text{and} \quad
\left\{ \begin{array}{l} y_4 = x_4 + x_5 + x_6, \\ y_5 = - x_4 + x_5, \\
y_6 = x_4,
\end{array} \right. $$ we have
$$ \begin{array}{ll} \sigma_1 : & y_2 \mapsto y_2 + y_1,\; y_3 \mapsto y_3 + y_2,\;
\text{and} \; y_1,y_4,y_5,y_6 \;\text{fixed,} \\ \sigma_2 : & y_5
\mapsto y_5 + y_4,\; y_6 \mapsto y_6 + y_5,\; \text{and} \;
y_1,y_2,y_3,y_4 \;\text{fixed,} \\ \,\tau \;: & y_1 \leftrightarrow
y_4,\; y_2 \leftrightarrow y_5,\; y_3
\leftrightarrow y_6, \\
\lambda^2 : & y_2 \mapsto - y_2,\; y_3 \mapsto y_3 + y_2,\; y_5
\mapsto - y_5,\; y_6 \mapsto y_6 + y_5,\; \text{and} \; y_1,y_4
\;\text{fixed.}
\end{array}
$$
Choose change of variables as follows: $$ \begin{array}{lll} z_1 =
y_1, & z_2 = y_2/y_1, & z_3 = y_3/y_1 + (y_2/y_1)^2 - y_2/y_1, \\
z_4 = y_4, & z_5 = y_5/y_4, & z_6 = y_6/y_4 + (y_5/y_4)^2 - y_5/y_4.
\end{array}
$$
We find $\sigma_1(z_2) = z_2 + 1$ and $\sigma_1$ fixes other
$z_i$'s, similarly $\sigma_2(z_5) = z_5 + 1$ and $\sigma_2$ fixes
other $z_i$'s, so by Artin-Schreier theory, $$ K(x_1,\dots,x_6)^{C_3
\times C_3} = K(z_1,\dots,z_6)^{C_3 \times C_3} = K(z_1,
z_2^3-z_2,z_3,z_4,z_5^3-z_5,z_6).
$$

Since $\tau : z_1 \leftrightarrow z_4, z_2 \leftrightarrow z_5, z_3
\leftrightarrow z_6$, we have $K(z_1,
z_2^3-z_2,z_3,z_4,z_5^3-z_5,z_6)^{\langle \tau \rangle} =
K(u_1,\dots,u_6)$, with $u_1 = z_1 + z_4, u_2 = (z_2^3-z_2) +
(z_5^3-z_5), u_3 = z_3 + z_6, u_4 = z_1 z_4, u_5 = z_1 (z_5^3-z_5) +
z_4 (z_2^3-z_2), u_6 = z_1 z_6 + z_4 z_3$.

One checks $\lambda^2 : u_2 \mapsto - u_2, u_5 \mapsto - u_5$, and
$u_1,u_3,u_4,u_6$ are fixed by $\lambda^2$, so finally $$
K(u_1,\dots,u_6)^{\langle \lambda^2 \rangle} =
K(u_1,u_2^2,u_3,u_4,u_2u_5,u_6),
$$ we are done.

\begin{remark} Over any field $K$ of characteristic $p > 0$, consider the
action of $C_p = \langle \sigma \rangle$ on $K(x_1,\dots,x_p)$ by
$\sigma : x_1 \mapsto x_2 \mapsto \cdots \mapsto x_p \mapsto x_1$
(as a generalization of the above action of $\sigma_1$ on
$x_1,x_2,x_3$). After making transformation $y_i = (-1)^{i-1}
\sum_{j=1}^{p} \binom{i+j-2}{i-1} x_j$ for $1 \le i \le p$ (taking
$\binom{k}{0} = 1$ for convention), the action becomes $\sigma(y_1)
= y_1$ and $\sigma(y_i) = y_i + y_{i-1}$ for $2 \le i \le p$. It can
be checked that $K(x_1,x_2,\dots,x_p)^{\langle \sigma \rangle} =
K(z_1,z_2,\dots,z_p)$, where (letting $\theta = y_2/y_1$) $z_1 =
y_1, z_2 = \theta^p-\theta$, and $z_k = y_k/y_0 - \binom{\theta}{1}
y_{k-1}/y_0 + \binom{\theta+1}{2} y_{k-2}/y_0 - \cdots + (-1)^{k-1}
\binom{\theta+k-2}{k-1}$ for $3 \le k \le p$. When $p = 3$, this
gives the result above for $K(x_1,x_2,x_3)^{\langle \sigma_1
\rangle}$. See \cite{Ku} for a more general situation.
\end{remark}

\section{Generic polynomials}

A positive answer to rationality problem for a finite group $G$ will
lead to a generic polynomial the group $G$. This polynomial
parametrises those Galois extensions with group $G$ in a suitable
sense, see \cite{Sa1}, \cite{De}, \cite{DeMc}, \cite{KM} for
details.

Let us consider our group $G_1 < S_6$ of order 72. To get a generic
polynomial for $G_1$, we simply write down the polynomial
\begin{multline*} (X - x_1)(X - x_2)(X - x_3)(X - x_4)(X - x_5)(X -
x_6) \\ = X^6 + a_1 X^5 + a_2 X^4 + a_3 X^3 + a_4 X^2 + a_5 X + a_6,
\end{multline*} and express those $a_i$'s explicitly as elements of
$K(x_1,\dots,x_6)^{G_1} = K(z_1,\dots,z_6)$. The result is as
follows, \begin{multline*} X^6 - z_1 X^5 + (z_2+z_4) X^4 - (z_3+z_5)
X^3 + \Big(z_6+\frac{z_4 z_2^2 - z_1 z_2 z_5 +
z_5^2}{4z_4-z_1^2}\Big) X^2 \\ - \frac{2z_4 z_2 z_3 - z_1 (z_2 z_6 +
z_5 z_3) + 2 z_5 z_6}{4z_4-z_1^2} X + \frac{z_4 z_3^2 - z_1 z_3 z_6
+ z_6^2}{4z_4-z_1^2} .
\end{multline*}

Now since $z_i$'s are homogeneous, we can give a simplified form by
letting $z_1 = 1$ (cf. \cite{KM}, Th. 7), and performing some change
of variables, the results are:

$\bullet$ If ${\rm char}(K) = 2$, a generic polynomial is
\begin{multline*} X^6 + X^5 + (t_1+t_3) X^4 + (t_2+t_4) X^3 + (t_5 +
t_3 t_1^2 + t_1 t_4 + t_4^2) X^2 \\ + (t_1 t_5 + t_4 t_2) X + (t_3
t_2^2 + t_2 t_5 + t_5^2) \;\;\in\; K(t_1,t_2,t_3,t_4,t_5)[X].
\end{multline*}

$\bullet$ If ${\rm char}(K) \ne 2$, a generic polynomial is
\begin{multline*} X^6 - 2 X^5 + (2 t_1 + t_3 + 1) X^4 - (2 t_2 + 2
t_4) X^3 + \Big(2 t_5 + t_1^2 +
\frac{(t_1 - t_4)^2}{t_3}\Big) X^2 \\
- \Big(t_1 t_2 + \frac{(t_1 - t_4) (t_2 - t_5)}{t_3}\Big) X +
\Big(t_2^2 + \frac{(t_2 - t_5)^2}{t_3}\Big) \;\;\in\;
K(t_1,t_2,t_3,t_4,t_5)[X].
\end{multline*}

By using Mathematica, we can get similar formulas for the other 3
groups $G_2,G_3,G_4$. But the results are too complicated, and
practically of little use.

\renewcommand{\refname}{\centering{References}}

\end{document}